%
\input amstex
\UseAMSsymbols
\input pictex
\vsize43pc \hsize28pc
\vsize=23truecm
\NoBlackBoxes
\parindent=18pt
  
   \font\rmk=cmr8    \font\itk=cmti8


\def\Hom{\operatorname{Hom}}

\def\Ext{\operatorname{Ext}}
\def\rad{\operatorname{rad}}
\def\add{\operatorname{add}}
\def\Ker{\operatorname{Ker}}
\def\Cok{\operatorname{Cok}}
\def\soc{\operatorname{soc}}
\def\Tr{\operatorname{Tr}}
\def\Im{\operatorname{Im}}

\def\arr#1#2{\arrow <1.5mm> [0.25,0.75] from #1 to #2}

\plainfootnote{}
{\rmk 2010 \itk Mathematics Subject Classification. \rmk 
Primary 
        16D90, 
        16G10. 
Secondary:
        16G70. 
}

\centerline{\bf Morphisms determined by objects:}
	\medskip
\centerline{\bf The case of modules over artin algebras.
}
                     \bigskip
\centerline{Claus Michael Ringel}     
		  \bigskip\medskip
{\narrower\narrower \noindent
Abstract. Let $\Lambda$ be an artin algebra. In his Philadelphia Notes,
M.~Auslander showed that any homomorphism between $\Lambda$-modules is 
right determined
by a $\Lambda$-module $C$, but a formula for $C$ which he wrote down has
to be modified. The paper includes complete and direct proofs of
the main results concerning right determiners of morphisms. We discuss the role
of indecomposable projective direct summands of a minimal right determiner
and provide a detailed  analysis of those morphisms which are right determined
by a module without any non-zero projective direct summand. This answers
a question raised in the book by Auslander, Reiten and Smal\o.
What we encounter is an intimate relationship to the vanishing of
$\Ext^2.$\par}
	\bigskip
Let $\Lambda$ be an artin algebra, the modules which we consider are
finitely generated left $\Lambda$-modules. A morphism
$\alpha\:X \to Y$ of $\Lambda$-modules is said to be {\it right
determined} by a $\Lambda$-module $C$ provided the following condition
is satisfied: given any morphism
$\alpha'\:X' \to Y$ such that $\alpha'\phi$ factors through $\alpha$ for
any $\phi\:C \to X'$, then $\alpha'$ itself factors through $\alpha.$
This definition is due to Auslander; the papers [A1] and [A2] 
are devoted to this concept. One of the main assertions of Auslander 
claims that any morphism
$\alpha\:X \to Y$  is right determined by $C = \Tr D(K) \oplus
P(Q)$, see [A2], Theorem 2.6;
here  $K$ is the kernel, $Q$ the cokernel of $\alpha$, and 
$\Tr (M)$ denotes the transpose, $D(M)$ the dual 
and $P(M)$ the projective cover of a module $M$. 

The aim of this note is to show that this assertion is not
correct as stated (in contrast to the weaker statements Theorem 3.17 (b) of
[A1] and Corollary XI.1.4 in [ARS]). 
In section 1, we will present corresponding examples. 
The assertion has to be slightly modified: not the projective cover of $Q$
is relevant, but the projective cover of the {\bf socle} $\soc Q$ of $Q$.
	\medskip

{\bf Theorem 1.} {\it Let $\alpha\:X \to Y$ be a morphism. Let $K$ be the kernel of
$\alpha$ and $Q$ the cokernel of $\alpha$. 
Then $\alpha$ is right determined by $\Tr D (K) \oplus P(\soc Q).$}
	\medskip
The modification of Auslander's treatment is formulated in Lemma 1 below (this
should replace [A2] Lemma 2.1.b).
Auslander's proof is somewhat hidden in two rather long papers, but there is a second
treatment of this topic in the book by Auslander, Reiten, Smal\o{} [ARS],
see the last chapter.
Still we feel that it may be appreciated if we provide a complete 
(and quite short) direct proof of Theorem 1. This will be done in section 2.
In section 3 we will use the same methods in order to describe the minimal right
determiner $T(\alpha)$ of $\alpha$, as it was introduced in [ARS].  
In section 4 we will discuss the following question: given a simple 
submodule $S$ of $\Cok(\alpha)$, when is $P(S)$ a direct summand 
of $T(\alpha)$ ? The final section 5 is devoted to a detailed analysis of the structure
of those maps $\alpha$ which are right determined by $\Tr D(K)$, with $K$ the kernel
of $\alpha,$ or, equivalently, by a module without an indecomposable projective
direct summand. The problem of characterizing this class was raised in [ARS].
		\medskip
{\bf Acknowledgment.}
Our interest in these questions was stimulated by a lecture of Henning Krause at the Shanghai Conference on Representation Theory of Algebras, October 2011, where he stressed
the relevance of Auslander's work, see also [K]. The author has to thank Hideto Asashiba for
having pointed out a wrong argument in a first version of the paper, as well
as Gordana Todorov for many helpful comments. 
	\bigskip\medskip
{\bf 1. Two Examples.}
	\medskip
{\bf Example 1.} 
Consider the quiver of type $\Bbb A_3$ with linear orientation, say
with simple modules indexed by $1,2,3$, such that $S(1)$ is projective,
$S(3)$ is injective. Let $\alpha\:S(1) \to P(3)$ be the inclusion map,
thus the kernel is zero, and the projective cover of the cokernel is again
$P(3)$. We claim that $\alpha$ is not right determined by $C = P(3).$
Consider the inclusion map $\alpha'\: P(2) \to P(3).$ Obviously, 
$\alpha'$ cannot be factored through $\alpha$. However, we have
$\Hom(C,P(2)) = 0,$ and the only map $\phi\:C \to P(2)$ (the zero-map)
has as composition with $\alpha'$  the zero-map
$C \to P(3)$. But the zero-map $C \to P(3)$ factors through $\alpha$, trivially.
	\medskip
{\bf Example 2.} 
Actually, an even easier example is given by the quiver $\Bbb A_2$, but here we
deal with $\alpha$ being a zero map (some may consider this as a degenerate case, thus
we presented first another example). Denote the two simple modules by
$S(1)$ and $S(2),$ with $S(1)$ being projective, $S(2)$ being injective.
We take as $\alpha$ the zero-map $0 \to P(2)$, its cokernel is $P(2)$ and already
projective. But $\alpha$ is not right determined by $C = P(2)$, since
the inclusion map $\alpha'\: S(1) \to P(2)$ does not factor through $\alpha$
(after all, $\alpha$ is zero), whereas 
for any map $\phi\:C \to S(1)$ (there is only the zero map) the composition
$\alpha'\phi$ factors through $\alpha.$ 
	\medskip
{\bf Remark.} Let us stress that Auslander's claim is correct in case 
$\Lambda$ is commutative, or, more generally, in case all the arrows of 
the quiver of $\Lambda$ are loops. Namely, in this case 
(and only in this case) $\add P(M) 
= \add P(\soc M)$ for any $\Lambda$-module $M$. 
	\bigskip\medskip
{\bf 2. The proof of Theorem 1.}
		\medskip
We start with the necessary amendment to Auslander's treatment.
	\medskip
Given an indecomposable projective module $P$, we always will denote the
inclusion map $\rad P \to P$ by $\iota$, the projection $P \to P/\rad P$ by $\pi$.
	\medskip
{\bf Lemma 1.} {\it Let $\alpha\: X \to Y$ be a morphism with image $\alpha(X)$.
Let $\alpha'\:X' \to Y$ be a morphism. 
Assume that for any simple submodule of the cokernel $Q = \Cok(\alpha)$ and any map $\phi\:P(S) 
\to X'$ with $\alpha'\phi(\rad P(S)) \subseteq \alpha(X)$, 
the map $\alpha'\phi$ factors through $\alpha$. Then the image of
$\alpha'$ is contained in $\alpha(X)$.}
	\medskip
Proof. We assume that the image of $\alpha'$ is not contained $\alpha(X)$ 
and want to derive a contradiction. 
Let us denote by $\gamma\:Y \to Q$ the cokernel map for $\alpha$.
By assumption, $\gamma\alpha' \neq 0$. Let $U$ be the image of $\gamma\alpha'$, with
epimorphism $\epsilon\:X' \to U$ and inclusion map $\mu\:U \to Q$, thus
$\mu\epsilon = \gamma\alpha'.$ Since $U$ is non-zero,
we may consider a simple submodule $S$ of $U$, say with inclusion map $\nu\:S \to U$.
Of course, $S$ is a simple submodule of $Q$. Let $\pi\:P(S) \to S$ be a projective cover of 
$S$. Since $P(S)$
is projective and $\epsilon$ is an epimorphism, we can lift $\nu\pi$
and obtain a map $\phi\:P(S) \to X'$ with $\epsilon\phi = \nu\pi$.
Note that 
$$
  \gamma\alpha'\phi = \mu\epsilon\phi = \mu\nu\pi.
$$ 
Since $\pi\iota = 0$, it follows that 
$$
   \gamma\alpha'\phi\iota = \mu\nu\pi\iota = 0.
$$
This shows that the image of $\alpha'\phi\iota$ is contained in the kernel
of $\gamma$, but this is $\alpha(X)$. In this way, we see that $\alpha'\phi(\rad P(S)) \subseteq \alpha(X).$

Thus, we are in the situation mentioned 
in the statement of the Lemma: 
there is given a map  $\phi\:P(S) \to X'$, such that $\alpha'\phi(\rad P(S)) \subseteq \alpha(X)$
and by the assumption  of the Lemma, we
know that the map $\alpha'\phi$ factors through $\alpha$, say $\alpha'\phi = 
\alpha\phi'$ for some $\phi'\:P(S) \to X.$ Therefore
$$
 \mu\nu\pi = \gamma\alpha'\phi = \gamma\alpha\phi' = 0,
$$
since $\gamma$ is the cokernel of $\alpha$. But $\mu\nu$ is a monomorphism, therefore
$\pi = 0,$ a contradiction. 
	\bigskip
Let us continue, as promised,  with the complete proof of Theorem 1.
The only prerequisite which we will
use is the existence of almost split sequences. To be precise: we will need for
any indecomposable non-injective module $M$ a non-split short exact sequence
$$
   0 @>>> M @>\sigma>> N @>\rho>> \Tr D(M) @>>> 0,
$$
such that for any map $\zeta\:M \to N'$ which is not a split monomorphism, there
is $\zeta'\:N \to N'$ with $\zeta = \zeta'\sigma.$ 
	\medskip
{\bf Lemma 2.} {\it Let $\alpha\: X \to Y$ be a morphism
with kernel $K$ and image $\alpha(X)$.
Let $\alpha'\:X' \to Y$ be a morphism with image contained in $\alpha(X)$. 
Assume that for any map $\phi\:\Tr D(K) 
\to X'$, the composition $\alpha'\phi$ factors through $\alpha$. 
Then $\alpha'$ factors through $\alpha$.}
	\medskip
{\bf Remark.} 
Given a morphism $\alpha\:X \to Y$, we may try to split off non-zero direct summands
of $X$ which lie in the kernel of $\alpha$. If this is not possible, then $\alpha$
is said to be {\it right minimal.} In general, we may write $X = X_0\oplus X_1$
with $X_0$ contained in the kernel of $\alpha$ and such that $\alpha|X_1$ is
right minimal; then we call the kernel of $\alpha|X_1$ the {\it intrinsic kernel}
of $\alpha$ (note that it is unique  up to isomorphism). An indecomposable direct
summand $L$ of the kernel of $\alpha$ is a direct summand of the intrinsic kernel,
if and only if the composition of the embeddings $L \subseteq K \subseteq X$ is not a split monomorphism.

It will be of interest in section 3 that one may replace in Lemma 2 the
kernel $K$ by the intrinsic kernel $K'$, thus the assertion of Lemma 2
can be strengthened as follows: 
{\it Assume that for any map $\phi\:\Tr D(K') 
\to X'$, the composition $\alpha'\phi$ factors through $\alpha$. 
Then $\alpha'$ factors through $\alpha$.}
	\medskip
Proof of Lemma 2 (and its strengthening). 
We may assume that $Y = \alpha(X)$, thus there is given the exact sequence $\eta$
with  epimorphism $\alpha\:X \to Y$ and kernel $\mu\:K \to X$.
We form the induced exact sequence $\eta'$ with respect to $\alpha',$ thus there
is the following commutative diagram with exact rows:
$$
\CD
 0 @>>> K @>\mu>> X @>\alpha>> Y @>>> 0  \qquad\qquad \eta \cr
   @.   @|       @AA\beta'A         @AA\alpha'A \cr
 0 @>>> K @>\nu>> W @>\beta>>      X' @>>> 0  \qquad\qquad \eta' 
\endCD
$$ 
If $\eta'$ is a split exact sequence, then $\alpha'$ factors through $\alpha$.

Let us assume that $\alpha'$ does not factor through $\alpha$, in order to derive a contradiction, again. Thus $\eta'$  is not a split exact sequence. Write $K = \bigoplus K_i$ with indecomposable
modules $K_i$ and projection maps $\pi_i\:K \to K_i.$ Since $\eta'$ does not split, there
is some index $i$ such that the exact sequence induced
from $\eta'$ by the map $\pi_i$ does not split. This means that we 
have the following commutative diagram
with exact rows which do not split:
$$
\CD
 0 @>>> K @>\nu>> W @>\beta>>      X' @>>> 0  \qquad\qquad \eta'  \cr
   @.   @V\pi_iVV       @V\pi'_iVV         @| \cr
 0 @>>> K_i @>\nu_i>> W_i @>\beta_i>>      X' @>>> 0  \qquad\qquad \eta'_i 
\endCD
$$

Let us add here, that $K_i$ has to be a direct summand of the intrinsic kernel of $\alpha$.
This observation is necessary in order to see that the remark made above is justified.

Since $\nu_i\:K_i \to W_i$ is a monomorphism which does not split,
we see that $K_i$ cannot be injective, thus there
is an almost split sequence
$$
 0 @>>> K_i @>\sigma_i>> V_i @>\rho_i>> \Tr D(K_i) @>>> 0,
$$ 
and $\nu_i$ can be factored as $\nu_i = \nu_i'\sigma_i$ for some $\nu_i'\:V_i \to W_i.$
Thus we obtain the following commutative square on the left, and therefore also
the map $\phi:\Tr D(K_i)\to X'$ with a commutative square on the right:
$$
\CD
 0 @>>> K_i @>\nu_i>> W_i @>\beta_i>>      X' @>>> 0  \qquad\qquad \eta'_i  \cr
   @.   @|       @AA\nu_i'A   @AA\phi A \cr
 0 @>>> K_i @>\sigma_i>> V_i @>\rho_i>>  \Tr D(K_i) @>>> 0   \qquad\qquad \omega_i
\endCD
$$
By assumption, the map $\phi\alpha'\: \Tr D(K_i) \to Y$ factors through $\alpha,$
that means there is $\phi'\: \Tr D(K_i) \to X$ with $\alpha\phi' = \alpha'\phi$.
Now, $W$ is the pullback of $\alpha, \alpha'$, thus there is a map $\phi''\:\Tr D(K_i) 
\to W$ such that $\beta\phi'' = \phi$ and $\beta'\phi'' = \phi'$. It follows that
$$
 \phi = \beta\phi'' = \beta_i\pi'_i\phi''.
$$
But if $\phi$ factors through $\beta_i$, then the exact sequence $\omega_i$
induced from $\eta'_i$ by $\phi$ has to split. This is a contradiction,
since $\omega_i$ is an Auslander-Reiten sequence, thus non-split.
	\bigskip
{\bf Proof of Theorem 1.} Let $\alpha\:X \to Y$ be a morphism with kernel $K$
and cokernel $Q$ and let $C = \Tr D (K) \oplus P(\soc Q).$
Let $\alpha'\:X' \to Y$ be a morphism such that $\alpha'\phi$
factors through $\alpha$ for any map $\phi\:C \to X'.$

If $S$ is a simple submodule of $Q$, then $P(S)$ 
is a direct summand of $P(\soc Q)$, thus of $C$. 
Thus, for any map $\phi\:P(S) 
\to X'$, the composition $\alpha'\phi$ factors through $\alpha$. Lemma 1 asserts that
the image of $\alpha'$ is contained in the image of $\alpha$.
Now we use that $\Tr D(K)$ is a direct summand of $C$,
thus for any map $\phi\:\Tr D(K)
\to X'$, the composition $\alpha'\phi$ factors through $\alpha$. 
Thus we can apply Lemma 2 in order to see that 
$\alpha'$ factors through $\alpha$. This shows that $\alpha$ is right determined by 
$C.$
	\bigskip
{\bf Example 3.} Let us add an example which may be illuminating, albeit it is extremely special.
Let $\Lambda$ be the path algebra of a finite directed quiver. Let $b$ be a vertex of the quiver and
assume that there are $s$ arrows starting in $b$, say $b \to a_i$ with $1\le i \le s$, and 
that there are $t$ arrows ending in $b$, say $c_j \to b$ with $1\le j \le t.$
For any vertex $x$, we denote by $S(x)$ the simple module with support $x$, by $P(x)$ 
the projective cover of $S(x)$, by $I(x)$ the injective envelope of $S(x)$. 

Let $\alpha$ be a non-zero map $X = P(b) \to I(b) = Y$, 
this is the homomorphism which we want to look at. Note that the image of $\alpha$ is $S(x).$
The kernel of $\alpha$ is the radical of $P(b)$, thus the direct sum of the modules
$P(a_i)$ with $1\le i \le s.$ The cokernel of $\alpha$ is the factor module of $I(b)$ modulo its socle, thus it is the direct sum of the modules $I(c_j)$ with $1\le j \le t.$ The projective
cover of the socle of $I(c_j)$ is $P(c_j)$. Altogether we see: the theorem 
asserts that $\alpha$ is right determined by the module 
$$
 C = \bigoplus_{i=1}^s \Tr D(P(a_i)) \oplus \bigoplus_{j=1}^t P(c_j).
$$ 
But this module $C$ is precisely the middle term of the almost split sequence starting in $P(b)$.

This should not come as a surprise. Namely, let $X'$ be an indecomposable module 
and assume that there is a non-zero map $\alpha'\:X' \to Y = I(b)$. Then there is a
map $\beta':P(b) \to X'$ with composition $\alpha'\beta' = \alpha.$ 
Now either $\beta'$ is invertible so that $\alpha'$ factors through $\alpha$,
or else $\beta'$ is not invertible and $\alpha'$ does not factor through $\alpha$.
In the latter case, 
$\beta'$ factors through the minimal left almost split map $\gamma\:P(b) \to C$ starting in $P(b)$, 
this means that there is some $\phi\:C \to X'$ with $\beta' = \phi\gamma$.
But if we look at the composition of $\phi$ and $\alpha'$, 
then one should be aware that no non-zero map $C \to I(b)$ factors through $\alpha$.
	\bigskip\medskip
{\bf 3. Minimal right determiners.}
	\medskip
Taking into account the Remark after Lemma 2, 
the Theorem we discuss can be strengthened as follows: {\it Any morphism
$X \to Y$ is right determined by $\Tr D(K')\oplus P(\soc Q)$,} where $K'$ is the
intrinsic kernel and $Q$ the cokernel of $\alpha.$ But one can do even better. 
	\medskip
Let us call a module $T = T(\alpha)$ a {\it minimal right determiner} for $\alpha$,
provided $T$ right determines $\alpha$ and is a direct summand of any module $C$
which right determines $\alpha$. 
According to [ARS], Proposition XI.2.4, a minimal right determiner for $\alpha$ 
exists and is the direct sum of all modules $N$ which almost factor through
$\alpha,$ one from each isomorphism class.
The aim of this section is to present a proof of this result using the considerations
of section 2.
	
We recall from [ARS] that an 
indecomposable module $N$ is said to {\it almost factor through $\alpha\:X \to Y$}
provided there is a morphism $\eta\:N \to Y$ which does not factor through $\alpha$
whereas for any radical map $\psi\:M \to N$, the composition $\eta \psi$ factors through 
$\alpha$. Obviously, the latter condition can be replaced by the condition that the map
$\eta\rho$ factors through $\alpha$, where $\rho$ is the minimal right almost split map ending
in $N$. Thus an indecomposable module $N$ almost factors though $\alpha$ provided
there exist a commutative diagram 
$$
\CD
   M @>\rho>> N \cr
   @V\eta'VV       @VV\eta V \cr
    X   @>\alpha>> Y
\endCD
$$
such that $\eta$ does not factor through $\alpha$ (with $\rho$  minimal right almost split).
Note that in case $N = P$ is 
(indecomposable) projective, the minimal right almost split map ending in $P$
is just the map $\iota\:\rad P \to P$.
	\medskip
{\bf Lemma 3.} {\it Let $P$ be an indecomposable projective module which almost factors through 
a map $\alpha$. Then $P$ is the projective cover of
a simple submodule of $\Cok(\alpha)$.}
	\medskip
Proof. 
Let $\eta\:P \to Y$ be a map which does not factor through
$\alpha$, whereas $\eta \iota$ factors through $\alpha$.
Consider the image $U$ of $\eta$ in $Y$ and the factor module $S = (U+\alpha(X))/\alpha(X)
\subseteq Y/\alpha(X) = \Cok(\alpha)$. Since $\eta(\rad P) \subseteq \alpha(X)$, we see
that $S$ is either simple or zero. But if $S = 0$, then $\eta(P) \subseteq \alpha(X)$
and the projectivity of $P$ implies that $\eta$ factors through $\alpha$. Since this is
not the case, $S$ is simple and $\eta$ provides an epimorphism $P \to S.$ 
	\bigskip
{\bf Lemma 4.} 
{\it Let $\alpha\:X \to Y$ be a morphism. Let $K'$ be the intrinsic kernel of
$\alpha$ and $P$ the direct sum of all indecomposable projective modules which almost factor
through $\alpha$, one from each isomorphism class. Then 
$\alpha$ is right determined by $\Tr D (K') \oplus P.$}
	\medskip
Proof: Let $\alpha'\:X' \to Y$ be a morphism which does not factor through $\alpha$. We have to
find an indecomposable module $C$ which is either of the form $\Tr D(L)$, where $L$ is 
a direct summand of $K'$ or a projective module
which almost factors through $\alpha$, and a morphism $\phi\:C \to X'$ such that $\alpha'\phi$
does not factor through $\alpha.$ According to the strengthened Lemma 2, such a pair
$C,\phi$ exists if the image of $\alpha'$ is contained in the image $\alpha(X)$ of $\alpha.$ 

Thus 
we can assume that the image of $\alpha'$ is not contained in $\alpha(X)$.
According to Lemma 1, there is a simple submodule $S$ of $Q$ and a map $\phi\:P(S) \to 
X'$ with $\alpha'\phi(\rad P(S)) \subseteq \alpha(X)$ such that $\alpha'\phi$ 
does not factor through $\alpha$. Write $\alpha = \alpha_2\alpha_1$ 
with inclusion map $\alpha_2\: \alpha(X) \to Y$. 
Using this notation,
$\alpha'\phi \iota = \phi'\alpha_2$ for some $\phi'$ (the restriction of $\phi$). 
If $\phi' = \alpha_1\phi''$, then $\iota\alpha'\phi = \alpha_2\alpha_1\phi'' = \alpha\phi''$
together with the fact that $\alpha'\phi$ does not factor through $\alpha$ shows that
$P(S)$ almost factors through $\alpha$, thus $P(S), \phi$ is the required pair.

Finally, we have to consider the case where $\phi'$ does not factor through $\alpha_1$. 
But then $\alpha_2\phi'$
does not factor through $\alpha$ (namely, $\alpha_2\phi' = \alpha\psi$ shows that
$\alpha_2\phi' = \alpha\psi = \alpha_2\alpha_1\psi,$ but $\alpha_2$ is injective, thus $\phi' = \alpha_1\psi$).
Now $\alpha_2\phi'$ is a morphism with image in $\alpha(X)$, thus as in the first part of the proof,
there is an indecomposable direct summand $C$ of  $K'$ and a map $\eta\:C \to \rad P$
such that $\alpha_2\phi'\eta$ does not factor through $\alpha$. If we rewrite
the composition $\alpha_2\phi'\eta = \alpha'\phi \iota \eta = \alpha'(\phi \iota\eta),$ then we 
see that we have achieved what we want, namely the pair $C, \phi \iota \eta$.  
	\bigskip
It remains to be seen that we have obtained in this may a minimal right determiner for
$\alpha$, at least up to multiplicities.
	\bigskip
{\bf Lemma 5.} {\it Assume that $\alpha$ is right determined by a module $C$. Let
$L$ be an indecomposable direct summand of the intrinsic kernel of $\alpha.$
Then $L$ is not injective, $\Tr D(L)$ is isomorphic to a direct summand of $C$,
and $\Tr D(L)$ almost factors through $\alpha$.}
	\medskip
Proof: 
Denote by $\mu\:K \to X$ and $\mu'\:L \to K$ the embeddings. And 
write $\alpha = \alpha_2\alpha_1$ with $\alpha_1\:X \to \alpha(X)$ surjective, and
$\alpha_2\:\alpha(X) \to Y$ the inclusion map. Since $\mu\mu'$ is an embedding which does not split, we see that $K$ is not injective,
thus there is an almost split sequence
$$
 0 @>>> L @>\sigma>> M @>\rho>> \Tr D(L) @>>> 0,
$$
and we can lift the map $\mu\mu'$ to $M$: there is a map $\mu''\:M \to X$ with $\mu''\sigma 
= \mu\mu'$. Since $\rho$ is the cokernel of $\sigma$, there is a map $\eta\:\Tr D(L) \to \alpha(X)$
such that $\eta \rho = \alpha\mu'',$ thus 
we obtain the following commutative diagram:
$$
\CD
 0 @>>> L @>\sigma>> M @>\rho>> \Tr D(L) @>>> 0 \cr
 @.    @V\mu'VV     @V\mu''VV   @VV\eta V     \cr
 0 @>>> K @>\mu>> X @>\alpha>> Y
\endCD.
$$ 

We claim that $\eta$ does not factor through $\alpha$. 
In order to proof this, we 
recall that $L$ is a direct summand of $K$, say $K = L\oplus L'$, and we  
form the induced exact sequence the given Auslander-Reiten sequence with the
split monomorphism $\mu'\:L \to K = L \oplus L'.$ The induced sequence is the
direct sum of the Auslander-Reiten sequence and a sequence of the form
$0 \to L'\to L'\to 0 \to 0$, in particular non-split, see the diagram below.
Since $\mu''\sigma = \mu\mu'$,
we obtain a map $\beta\:M\oplus L'\to X$ and then a map $\beta'\:\Tr D(L) \to \alpha(X)$
such that the following diagram is commutative:
$$
\CD
 0 @>>> L @>\sigma>> M @>\rho>> \Tr D(L) @>>> 0 \cr
 @.    @V\mu'VV     @V\mu''VV   @|     \cr
 0 @>>> K @>\mu>> M\oplus L' @>>> Y @>>> 0 \cr
 @.    @|         @V\beta VV   @VV\beta' V     \cr
 0 @>>> K @>\mu>> X @>\alpha_1>> \alpha(X) @>>> 0
\endCD
$$
Note that a comparison with the diagram above shows that $\eta = \alpha_2\beta'.$
From the diagram we see that the horizontal middle sequence is induced from the
lower sequence by $\beta'$. Since the horizontal middle sequence does not split,
we see that $\beta'$ does not factor through $\alpha_1$.
Now assume that $\eta$ factors through $\alpha,$ say $\eta = \alpha\zeta$
for some $\zeta\:\Tr D(L) \to X$. Then
$$
 \alpha_2\alpha_1\zeta = \alpha\zeta = \eta = \alpha_2\beta',
$$
implies that $\alpha_1\zeta = \beta'$, since $\alpha_2$ is injective.
But we know already that $\beta'$ does not factor through $\alpha_1$,
thus $\eta$ does not factor through $\alpha,$ as we wanted to show.

Since $C$ right determines $\alpha,$ and $\eta\:\Tr D(L) \to Y$ does not factor through
$\alpha$, there has to exist a morphism
$\phi\:C \to \Tr D(L)$ such that also $\eta\phi$ cannot be
factored through $\alpha$. Now again we use that the upper
sequence is an Auslander-Reiten sequence. Assume that $\phi$ is
not split epi. Then there is $\phi'\:C \to M$ such that $\rho\phi' = \phi$,
and therefore
$$  
 \eta\phi = \eta\rho\phi' = \alpha\beta\mu''\phi'
$$
is a factorization of $\eta\phi$ through $\alpha,$ a contradiction.
This shows that $\phi$ is split epi, thus $\Tr D(L)$ is isomorphic to 
a direct summand of $C$.

Finally, we see that $\Tr D(L)$ almost factors through $\alpha$, since there is
the diagram
$$
\CD
 M @>\rho>> \Tr D(L) \cr
 @V\mu''VV         @VV\eta V  \cr
 x @>\alpha>>       Y 
\endCD
$$
and $\eta$ does not factor through $\alpha.$
	\bigskip
{\bf Lemma 6.} {\it Assume that $\alpha$ is right determined by a module $C$. Let
$P$ be an indecomposable projective which almost factors through $\alpha$
Then $P$ is isomorphic to a direct summand of $C$.}
	\medskip
Proof: There exists a commutative diagram
$$
\CD
  \rad P @>\iota>> P \cr
   @VVV       @VV\eta V \cr
    X   @>\alpha>> Y
\endCD
$$
such that $\eta$ does not factor through $\alpha.$ 
Since $C$ right determines $\alpha$,
there must exist $\phi\:C \to P$ such that also $\eta\phi$ does not factor through $\alpha$.
Now $\phi$ does not map into $\rad P$, since otherwise $\eta\phi$ would factor through $\alpha.$
But this means that $\phi$ is surjective and therefore a split epimorphism.
	\bigskip
{\bf Theorem 2.} {\it Let $\alpha\:X \to Y$ be given. Let
$T$ be the direct sum of modules of the form $\Tr D(L)$, where $L$
is an indecomposable direct summand of the intrinsic kernel of
$\alpha$ and of the 
indecomposable projective modules which almost factor
through $\alpha$, one from each isomorphism class. Then 
$T$ is a minimal right determiner for $\alpha.$}
	\bigskip
Proof. This is a direct consequence of the Lemmata 4, 5 and 6.
	\medskip
{\bf Corollary 1.} {\it Let $\alpha\:X \to Y$ be given. 
A non-projective indecomposable module $N$ almost factors through $\alpha$
if and only if $N = \Tr D(L)$ for some indecomposable direct summand $L$ 
of the intrinsic kernel of $\alpha$.}
	\medskip
Proof. On the one hand, we have seen in Lemma 5 that the modules of the form $\Tr D(L)$
almost factor through $\alpha$. On the other hand, it is clear that 
an indecomposable module which almost factors through 
$\alpha$ is a direct summand of any right determiner for $\alpha$
(see for example [ARS] Lemma XI.2.1), thus of $T(\alpha)$.
	\bigskip
{\bf  Corollary 2.} {\it Let $\alpha\:X \to Y$ be given. 
An indecomposable module $N$ almost factors through $\alpha$ if and
only if it is a direct summand of $T(\alpha).$}
	\bigskip

	\bigskip\bigskip
{\bf 4. The indecomposable projective direct summands of $T(\alpha)$.}
	\medskip
Theorem 2 shows that $T(\alpha)$ has two kinds of indecomposable direct summands:
First of all, there are those 
of the form $\Tr D(L)$, where $L$ is any direct summand of the intrinsic
kernel of $\alpha$, and clearly they are never projective. 
Second, there may be indecomposable projective modules.
Here we want to discuss these latter summands.
	\bigskip
Recall that if $S$ is a simple module such that $P(S)$ 
is a direct summand of $T(\alpha)$,
then, according to Lemma 3, $S$ is a simple
submodule of $\Cok(\alpha).$ But the converse does not hold.
{\it Not every module $P(S)$ with $S$ a simple
submodule of $\Cok(\alpha)$ almost factors through $\alpha$.}
	\medskip
{\bf Example 4.} This example has been exhibited in the book
of Auslander, Reiten, Smal\o{} [ARS], after Proposition XI.1.6.
Let $\Lambda$ be a local uniserial ring with the unique simple module
$S$, and let $\alpha\:P \to Y$ be a morphism with $P$ the 
indecomposable projective module and
$Y$ also indecomposable. If $P = P(S)$ almost factors through $\alpha$,
then $\alpha = 0$, and therefore $\alpha$ is right determined by 
$\Tr D(\Ker(\alpha)).$ 
	\medskip
Actually, for any artin algebra with global dimension at least 2 
there do exist corresponding examples, as the following basic observation
shows:
	\medskip
{\bf Example 5.} {\it Let $\delta\:P_1 \to P_0$ be a minimal presentation
of a simple module $S$. If $P(S) (= P_0)$ almost factors through $\delta$,
then $\delta$ is injective,} thus the projective dimension of $S$ is
at most 1. Proof: Write $\delta = \iota\epsilon$, where $\iota\:\rad P_0 \to P_9$
is the inclusion map. If $P_0$ almost factors through $\delta$, there is
$\eta\:P_0 \to P_0$ not factoring through $\delta$ and $\eta'\:\rad P_0 \to P_1$
such that $\eta\iota = \delta\eta'$, whereas $\eta$ does not factor through $\delta$.
Then $\delta$ does not map into $\rad P_0$, therefore $\eta$ has to be invertible,
and $\eta\iota = \iota\epsilon\eta'$ implies that $\iota = \eta^{-1}\iota\epsilon\eta'$,
thus $1_{\rad P_0} = \epsilon\eta'$. But  this means that $\epsilon$ is  split
epimorphism, thus an isomorphism (since it is a projective cover).
	\bigskip
Here are three sufficient conditions for $P(S)$ to be a direct
summand of $T(\alpha)$.
	\medskip
{\bf Proposition 1.} {\it Let $\alpha\: X \to Y$ be a monomorphism with cokernel $Q$.
If $S$ is a simple submodule of $Q,$
then $P(S)$ almost factors through $\alpha$.}
	\medskip
Proof. We may assume that $\alpha$ is an inclusion map.
Since $S$ is a submodule of $Y/X$, there is a map $\eta\:P(S) \to Y,$
such that the composition of $\eta$ with $Y \to Y/X$ maps onto $S$.
But then $\eta(\rad P(S)) \subseteq X.$ Thus $P(S)$ almost
factors through $\alpha$.
	\medskip
{\bf Proposition 2.} {\it Let $\alpha\: X \to Y$ be a morphism.
If $S$ is a simple submodule of $Y$ with $S\cap \alpha(X) = 0,$
then $P(S)$ almost factors through $\alpha$.}
	\medskip
Proof: Let $S$ is a simple submodule of $Y$, and  
let $\eta\:P(S) \to Y$ be a morphism with image $S$. Then $\eta\iota = 0$.
Thus the following
diagram commutes:
$$
\CD
  \rad P(S) @>\iota>> P(S) \cr
   @V0VV       @VV\eta V \cr
    X   @>\alpha>> Y 
\endCD .
$$
Since $S\cap \alpha(X) = 0$, we see that $\eta$ does not factor through $\alpha.$
	\medskip
{\bf Proposition 3.} {\it Let $\alpha\:X \to Y$ be a morphism with cokernel $Q$.
Let $S$ be a simple submodule of $Q$. If the projective dimension of $S$ is 
at most $1$, then $P(S)$ almost factors through $\alpha$.}
	\medskip
{\bf Proof.}
Let $\pi\:P(S) \to S$ be a projective cover and $\nu\:S\to Q$ 
the inclusion map. Let $\gamma\:Y \to Q$ be the cokernel map. 
The projectivity of $P(S)$ yields a map $\eta\:P(S) \to Y$
such that 
$\gamma\eta = \nu\pi.$ 
We denote the projection $Y \to Y/\alpha(X) = Q$ by $\gamma$. Then
$\gamma\eta\iota = \nu\pi\iota = 0,$ thus $\eta$ maps $\rad P(S)$ into $\alpha(X).$
This shows that we have the following commutative diagram 
$$
\CD
  \rad P(S) @>\iota>> P(S) \cr
   @V\eta'VV       @VV\eta V \cr
    \alpha(X)  @>\alpha_2 >> Y 
\endCD ,
$$
as before we write $\alpha = \alpha_2\alpha_1$ where $\alpha_2\:\alpha(X) \to Y$ 
is the canonical inclusion of $\alpha(X) = \alpha(X)$ into $Y$.
Since the projective dimension of $S$ is at most $1$,
we know that $\rad P(S)$ is projective, thus we can lift $\eta'$ and obtain
$\eta''\:\rad P(S)\to X$ with $\alpha_1\eta'' = \eta'$, thus there is the commutative diagram
$$
\CD
  \rad P(S) @>\iota>> P(S) \cr
   @V\eta''VV       @VV\eta V \cr
    X   @>\alpha>> Y 
\endCD .
$$
Of course, $\eta$ does not factor through $\alpha$ since $\gamma\eta \neq 0.$
	\medskip
It follows that for a hereditary artin algebra, the projective cover $P(S)$ of any
simple submodule of $\Cok(\alpha)$ is a direct summand of $T(\alpha).$
	\medskip
Finally, there is the following characterization:
	\medskip
{\bf Proposition 4.} {\it Let $S$ be a simple module. Then 
$P(S)$ is a direct summand of $T(\alpha)$ if and only if there exists
a module $J$ with submodule $X$ and $J/X = S$ and
a morphism $\widetilde\alpha\:J \to Y$ such that its restriction to
$X$ is $\alpha$ and the kernels of $\alpha$ and $\widetilde \alpha$ coincide.} 
The condition that the kernels of $\alpha$ and $\widetilde \alpha$ coincide
is equivalent to the condition that the image of $\alpha$ is properly
contained in the image of $\widetilde \alpha.$
	\medskip
Proof: First, let us assume that there exists 
a module $J$ with submodule $X$ and $J/X = S$ and
a morphism $\widetilde\alpha\:J \to Y$ such that its restriction to
$X$ is $\alpha$ and such that the image of $\alpha$ is a properly
contained in the image of $\widetilde \alpha.$
Denote the projection map $J \to J/X = S$ by $\epsilon.$ 
Let $\pi\:P(S) \to S$ be a projective cover and lift it to $J,$ thus we obtain
$\pi'\:P(S) \to J$ such that $\epsilon\pi' = \pi.$
Since $\epsilon\pi'(\rad P(S)) = \pi(\rad P(S)) = 0$, we have
$\pi'(\rad P(S)) \subseteq X$. Let us denote by $\pi''\:\rad P(S) \to X$
the restriction of $\pi'$ to $\rad P(S)$. Then the diagram
$$
\CD
  \rad P(S) @>\iota>> P(S) \cr
   @V\pi''VV       @VV\widetilde\alpha\pi' V \cr
    X  @>\alpha >> Y 
\endCD 
$$
commutes, since $\widetilde\alpha|X = \alpha.$

It remains to be seen that $\widetilde\alpha\pi'$ does not
factors through $\alpha$. Assume for the contrary that $\widetilde\alpha\pi' = \alpha\zeta$,
for some map $\zeta\:P(S) \to Y.$ 
Now $J = X + \pi'(P(S)),$ thus
$$
\align
 \widetilde\alpha(J) &= \widetilde\alpha(X + \pi'(P(S))) =
 \alpha(X) + \widetilde\alpha\pi'(P(S)) \cr
  &=
 \alpha(X) + \alpha\zeta(P(S)) = \alpha(X),
\endalign
$$
contrary to our assumption. 

Conversely, assume that $P(S)$ almost factors through $\alpha$, thus we have 
a diagram of the following form
$$
\CD
  \rad P(S) @>\iota>> P(S) \cr
   @V\eta'VV       @VV\eta V \cr
   X  @>\alpha >> Y 
\endCD
$$  
and $\eta$ does not factor through $\alpha$, thus the image of $\eta$ is not
contained in the image of $\alpha$.
Starting with the exact sequence with monomorphism $\iota$, we form the
sequence induced by $\eta'$ and obtain the following commutative diagram with
exact rows:
$$
\CD
0 @>>>   \rad P(S) @>\iota>> P(S) @>\pi>> S @>>> 0 \cr
@.        @V\eta' VV               @VV\eta''V         @|  \cr
0 @>>>       X  @>\iota' >> J   @>>>    S @>>> 0
\endCD
$$
Since $\eta\iota = \alpha\eta'$, there is a map $\widetilde\alpha\:J \to Y$
such that $\alpha = \widetilde\alpha\iota'$ and $\eta =\widetilde\alpha\eta''$. 
Thus, we see that
$\alpha$ has an extension $\widetilde\alpha$ to $J.$ 
Since $\eta =\widetilde\alpha\eta''$, the image of $\eta$ is contained in the
image of $\widetilde\alpha$. This shows that 
the image of $\widetilde\alpha$ cannot be equal to the image of $\alpha$,
since otherwise the image of $\eta$ would be contained in the
image of $\alpha$, in contrast to our assumption.
This concludes the proof.
	\bigskip
Proposition 4 (but also already Proposition 3) show that
the obstructions for the projective cover $P(S)$ of a simple submodule of $\Cok(\alpha)$
to be a direct summand of $T(\alpha)$ are elements of $\Ext^2$, namely the 
equivalence classes of the exact sequences
$$
 0 @>>> K @>>> X @>>> J @>>> S @>>> 0, \tag{$*$}
$$
where $K$ is the kernel of $\alpha$ and 
$J = \gamma^{-1}(S)$ (here  $\gamma$ is the cokernel map $Y \to \Cok(\alpha)$)
and where the composition of the map $X \to J$ with the inclusion map $J
\to Y$ is just $\alpha$. Thus we have:
	\medskip
{\bf Corollary.} {\it Let $\alpha\:X \to Y$ be a morphism with kernel $K$ and cokernel $Q$. 
If $S$ is a submodule of $Q$ and $\Ext^2(S,K) = 0$, then $P(S)$ is a direct summand of $T(\alpha).$}
	\bigskip\bigskip
{\bf 5. Kernel-determined morphisms.}
	\medskip
Since any morphism $\alpha$ is right determined by the direct sum of the module
$\Tr D(\Ker(\alpha))$
and a projective module $P$, one may ask for a characterization of those morphisms $\alpha$
for which one of these two summands already right determines $\alpha$.
	\medskip
First, let us deal with the morphisms which are right determined by a projective module.
Here, the answer is well-known and easy to obtain:
{\it A morphism $\alpha$ is right determined
by a projective module if and only if $\alpha$ is injective} (see Theorem 1 and Lemma 5). 

Also, {\it an inclusion map
$X \to Y$ is right determined by the projective module $P$, if and only if $P$ 
generates the socle of $Y/X$.} (If $P$ generates the socle of $Y/X$, 
then $P$ right determines $\alpha$
according to Theorem 1. Conversely, assume that $P$ right determines $\alpha$, and 
let $S$ be a simple submodule of $Y/X$. According to Proposition 1,
$P(S)$ almost factors through $\alpha,$ thus Theorem 2 asserts that $P(S)$ is a direct
summand of $P$. This shows that $P$ generates the socle of $Y/X$.)
There is the following consequence: 
{\it If we fix a projective module $P\neq 0$, and consider any module $X$, then
there are morphisms $\alpha\:X \to Y$ with $Y$ of arbitrarily large length, such that
$\alpha$ is right determined by $P$} (just take the inclusion maps of the form $X \to Y$
with $Y$ the direct sum of $X$ and arbitrarily many copies 
of $P/\rad(P)$). If $\Lambda$ is
representation-infinite, then there are even such examples with $Y$ indecomposable.
	\medskip
The second case are the morphisms $\alpha$ which are right determined by $\Tr D(\Ker(\alpha))$,
we call them {\it kernel-determined} morphisms.
This is the topic of the considerations in this section. Note that the problem of
characterizing these maps has been raised already in [ARS], 368-369. 
	\medskip
{\bf Lemma 7.} {\it Let $\alpha$ be a morphism.
The following conditions are equivalent:
\item{\rm (i)} $\alpha$ is right determined by $\Tr D (K)$, where $K$ is the kernel of $\alpha$.
\item{\rm (ii)} $\alpha$ is right determined by $\Tr D (K')$, where $K'$ is the intrinsic
  kernel of $\alpha$.
\item{\rm (iii)} $\alpha$ is right determined by a module $C$ without an indecomposable projective
direct summand.\par}
	\medskip
Proof. Clearly (ii) $\implies$ (i) $\implies$ (iii). Now assume (iii). According to Theorem 2,
any indecomposable projective module $P$ which almost factors through $\alpha$ is a direct summand
of $C$, thus there are no such modules $P$. Using again Theorem 2, we see that (ii) is
satisfied.
	\medskip
Note that $\alpha$ is kernel-determined if and only if 
the equivalent conditions of lemma 7 are satisfied.
Let us first show that for a kernel-determined morphism $\alpha\:X \to Y$,
the length of $Y$ is bounded by a number which only depends on $X$.
We denote by $|M|$ the length of the module $M$.
	\medskip
{\bf Lemma 8.} 
{\it If $\alpha\: X \to Y$ is kernel-determined, then $Y$ is an
essential extensions of $\alpha(X)$; 
in particular, $|Y| \le q|X|$ where
$q$ is the maximal length of an indecomposable injective module.}
	\medskip
If $Y$ is an essential extension of $N = \alpha(X)$, then we may assume that
$Y$ is a submodule of $I(N)$ with $N \subseteq Y.$ 
	\medskip
Proof of lemma 8. According to Proposition 2, there is no simple submodule $S$ of $Y$
with $S\cap \alpha(S) = 0,$ this jut means that $Y$ is an essential extension of $\alpha(X)$.
Thus $Y$ can be considered as a submodule of the injective envelope $I$ of $\alpha(X)$.
But then $|I| \le q|\alpha(X)| \le q|X|.$
	\medskip 
Given a module $M$, let $\overline M$ be a module having $M$ as an essential submodule 
with $\overline M/M$ semisimple and such that $\overline M$ is of maximal possible length;
we call $\overline M$ a {\it small envelope} of $M$.
We can construct $\overline M$ as follows:
$$
  \overline M  =\omega^{-1}(\soc I(M)/M),
$$ 
where $I(M)$ is an
injective envelope of $M$ and $\omega\:I(M) \to I(M)/M$ is the canonical projection map
(thus, if necessary, we will assume that $\overline M$ is a submodule of $I(M)$ which contains
$M$).
Clearly, any homomorphism $\phi\:M \to N$ gives rise to an extension $\overline\phi\:
\overline M \to \overline N$ (by this we mean a homomorphism 
whose restriction to $M$ is just $\phi$). 
Let us stress that usually $\overline\phi$ is not uniquely determined 
(the construction $M \mapsto \overline M$ is not functorial). But there is the following
unicity result which is of interest for the further considerations:
	\medskip
{\bf Lemma 9.} {\it Let $\epsilon\:X \to N$ be an epimorphism, and choose 
an injective envelope $I(N)$ of $N$. Then there is a (uniquely determined)
submodule $N \subseteq I_\epsilon(N) 
\subseteq I(N)$ with the following property: If $\overline X$ is a small envelope of $X$
and $\overline \epsilon\: \overline X \to I(N)$ is an extension of $\epsilon$, then
$\overline\epsilon(\overline X) = I_\epsilon(N).$}
	\medskip
Proof: If we deal with two extensions of $\epsilon,$ say $\epsilon_1, \epsilon_2: 
\overline X \to I(N)$, then the difference
$\epsilon_2-\epsilon_1$ vanishes on $X$ and its image is a semisimple module. But any semisimple
submodule of $I(N)$ is contained in $N$ and $N = 
\epsilon(X) \subseteq \epsilon_1(\overline X)$.
Thus, $\epsilon_2 = \epsilon_1 + (\epsilon_2-\epsilon_1)$ 
shows that
$$
 \epsilon_2(\overline X) \subseteq \epsilon_1(\overline X) + (\epsilon_2-\epsilon_1)(\overline X) 
   \subseteq  \epsilon_1(\overline X) + N \subseteq \epsilon_1(\overline X).
$$
Of course, by symmetry we also have $\epsilon_2(\overline X) \subseteq \epsilon_1(\overline X),$
and therefore equality.
	\medskip
Clearly, the submodule $I_\epsilon(N)$ incorporates the information about the
vanishing in $\Ext^2$ of the exact sequences of the form $(*)$, where $K\to X$ is the kernel
map for $\epsilon\:X \to N$.
	\bigskip
{\bf Theorem 3.} {\it Let $\epsilon\: X \to N$ be an epimorphism. 
Consider a submodule $N \subseteq Y \subseteq I(N)$ and denote by
$\nu\:N \to Y$ the inclusion map. Let $\alpha = \nu\epsilon$. Then 
$\alpha\:X \to Y$ is kernel-determined
if and only if $Y \cap I_\epsilon(N) = N.$}
	\medskip
Proof. We fix some notation. Let $D = Y \cap I_\epsilon(N)$.
Let $\nu'\:N \to D$, $\nu''\:D \subseteq Y$, $\nu'''\:Y \to I(N),$
$\kappa\:D \to I_\epsilon(N)$, and 
$\mu\:X \to \overline X$ be the inclusion maps. Thus we have $\nu = \nu''\nu$.

The inclusion map $\kappa\nu'\:X \to I_\epsilon(N)$ is part of the
following commutativity relation:
$$
\kappa\nu'\epsilon = \overline\epsilon_1\mu, 
 \tag{1}
$$ 
where we denote by $\overline\epsilon_1$ the epimorphism part of an extension 
$\overline\epsilon$ of $\epsilon.$
	
First, let us assume that $\nu'\: N\subset D = Y \cap I_\epsilon(N)$
is a proper inclusion. Then there exists an indecomposable
projective module $P$ and a homomorphism $\eta\:P\to D$ such that the image of $\eta$ 
does not lie inside $N$. Now $\overline\epsilon_1\:X \to  I_\epsilon(N)
$ is surjective, thus we can lift the map $\kappa\eta\:P(S) \to I_\epsilon(N)$
to $\overline X$ and obtain $\eta'\:P(S) \to \overline X$ such that
$$
 \overline\epsilon_1\eta' = \kappa\eta \tag{2}
$$

Also note that $\eta'\iota$ maps into the radical of $\overline X$, thus into $X$.
This shows that there is $\eta''\:\rad P(S) \to X$ such that
$$
 \mu\eta'' = \eta'\iota. \tag{3}
$$
Altogether, we deal with the following diagram:
$$
{\beginpicture
\setcoordinatesystem units <2cm,1.5cm>
\put{$\rad P(S)$} at 0 2
\put{$P(S)$} at 1 2
\put{$D$} at 2 2
\put{$X$} at 0 1
\put{$\overline X$} at 1 1
\put{$N$} at 0 0
\put{$I_\epsilon(N)$} at  2 0
\arr{0.4 2}{0.8 2}
\arr{1.2 2}{1.8 2}
\arr{0.2 1}{0.8 1}
\arr{0.2 0}{1.7 0}
\arr{2 1.8}{2 0.2}
\arr{0 0.8}{0 0.2}
\arr{1.2 .8}{1.8 0.2}
\arr{0 1.3}{0 1.2}
\arr{1 1.3}{1 1.2}
\setdashes <1mm>
\plot 0 1.8  0 1.3 /
\plot 1 1.8  1 1.3  /
\put{$\iota\strut$} at 0.5 2.2 
\put{$\eta''\strut$} at -.15 1.5 
\put{$\epsilon\strut$} at -.15 .5 
\put{$\eta'\strut$} at 1.15 1.5 
\put{$\kappa\strut$} at 2.15 1 
\put{$\kappa\nu'\strut$} at 1 -.2 
\put{$\overline \epsilon_1\strut$} at 1.6 0.6 
\put{$\mu\strut$} at .5 1.15 
\put{$\eta\strut$} at 1.5 2.2 
\put{$(1)$} at 0.7 0.45
\put{$(2)$} at 1.55 1.3
\put{$(3)$} at 0.5 1.55
\endpicture}
$$

Using the three equalities (1), (3), (2), we see:
$$
 \kappa\nu'\epsilon\eta'' = \overline\epsilon_1\mu\eta'' = \overline\epsilon_1\eta'\iota = 
  \kappa\eta\iota.
$$
but $\kappa$ is injective, thus 
$\nu'\epsilon\eta''  =  \eta\iota,$ and therefore
$$
 \alpha\eta'' = \nu''\nu'\epsilon\eta''  =   \nu''\eta\iota.
$$
This asserts that the following diagram commutes
$$
\CD
  \rad P @>\iota>> P \cr
   @V\eta''VV           @VV\nu''\eta V \cr
    X @>\alpha>> Y 
\endCD
$$
Since by construction the right map $\nu''\eta$ does not map into $N$, it does
not factor through  $\epsilon$, thus also not through 
$\alpha = \nu''\nu'\epsilon$, therefore we see that $P$ almost factors through $\alpha$.
But this shows that $\alpha$ is not kernel-determined.

Conversely, let us assume that $\alpha = \nu''\nu\epsilon$ is not 
kernel-determined, thus there is an indecomposable projective module $P$ and a
commutative diagram
$$
\CD
 \rad P @>\iota>> P \cr
  @V\psi'VV      @VV\psi V \cr
   X @>\alpha>> Y
\endCD
$$
such that $\psi$ does not factor through $\alpha = \nu\epsilon$, thus $\psi(P)$
is not contained in $N$. Let us form a pushout of $\iota$ and $\psi'$, say
$$
\CD
 \rad P @>\iota>> P \cr
  @V\psi'VV      @VV\psi'' V \cr
   X @>\iota'>> J
\endCD,
$$
we obtain a map $\beta\:J \to Y$ such that $\beta\psi'' = \psi$
and $\beta\iota' = \alpha.$ Since $Y$ is a submodule of $I(N)$, the image
$\beta(J)$ of $\beta$ is a submodule of $Y$, thus of $I(N)$.

Let us show that $\iota'$ does not split and its cokernel is simple.
The cokernel of $\iota'$ is isomorphic to the cokernel of $\iota$,
thus simple. 

Let us show that {\it the kernel of $\beta$ is just $\iota(K)$,} where $K$ is the kernel of 
$\alpha$. Since $\beta\iota = \alpha$, we see that $\iota(K)$ is contained in
the kernel of $\beta$, thus it remains to show that $|\Ker(\beta)| \le |\Ker(\alpha)|$
(note that $\iota$ is injective). 
Since $\alpha = \beta\iota$, the image $N$ of $\alpha$ is contained in the image of $\beta$.
This must be a proper inclusion. Otherwise, we use $\psi = \beta\psi''$ in order to obtain
that $\Im(\psi) \subseteq \Im(\beta) = \Im(\alpha) = N$, a contradiction. 
Thus $|\Im(\beta)| \ge |\Im(\alpha)|+1.$    
Therefore
$$
 |\Ker(\beta)| = |J|-|\Im(\beta)| \le |X|+1 - |\Im(\alpha)| -1 = |\Ker(\alpha)|.
$$

It follows that  $\iota'$ does not split. Otherwise
we have $J = \iota(X)\oplus S$. Now the kernel of $\beta$ is $\iota(K)
= \iota(K) \oplus 0$, and therefore $\beta$ would provide an embedding of 
$X/K \oplus S$ into $Y$. However, by assumption, $Y$ is an essential extension of $N = X/K$,
a contradiction.

Thus we have shown that $\iota'$ is a monomorphism with simple cokernel, and it does not split.
Therefore, we may assume that $J$ is a submodule of $\overline X$,
If we compose $\beta$ with $\nu'''\:Y \to I(N)$, we obtain the following
commutative square
$$
\CD
 X @>\iota'>> J \cr
 @V\epsilon VV      @VV\nu'''\beta V \cr
 N @>\nu'''\nu>> I(N)
\endCD
$$
which shows that $\nu'''\beta$ is part of an extension $\overline\epsilon\:\overline X \to
I(N)$ of $\epsilon.$ As a consequence, the image of $\beta$ is contained in $I_\epsilon(N).$
But the image
$\beta(J)$ of $\beta$ is also a submodule of $Y$, that $\beta(J)
\subseteq D.$ 

Since $\beta\psi'' = \psi$, the image of $\beta$ contains the image of $\psi$,
thus $\beta(J)$ is not contained in $N$. 

Altogether we see that $\beta(J) \subseteq Y\cap I_\epsilon(N)$,
and $\beta(J) \not\subseteq N$, thus $Y\cap I_\epsilon(N) \neq N.$
This completes the proof.
	\bigskip
{\bf Example 4, continued.} Again, let $\Lambda$ be a local uniserial
ring.
Let $X, Y$ be indecomposable $\Lambda$-modules and $\alpha\:X \to
Y$ a morphism. We have noted above that if $X$ is projective and $\alpha\neq 0$, 
then $\alpha$ is kernel-determined. On the other hand, if $\alpha$ is surjective,
then again, $\alpha$ is kernel-determined. 

But also the converse is true: If $\alpha\:X \to Y$ is kernel-determined, then
either $\alpha \neq 0$ and $X$ is projective, or else $\alpha$ is surjective.
Here is the proof: Assume that $\alpha$ is kernel-determined. 
According to Proposition 2, we must have $\alpha\neq 0.$  
Assume that $X$ is not projective, thus also not injective. 
Write $\alpha = \nu\epsilon$, where $\epsilon$ is surjective and 
$\nu\:N \to Y$ is the inclusion of a non-zero submodule $N$ of $Y$. 
Since $X$ is not injective, $X$ is a proper submodule of $\overline N$.
Let $\overline\epsilon\:\overline X \to \overline N$ be an extension of $\epsilon$.
Then also $\overline\epsilon$ is surjective. But this means that
$I_\epsilon(N) = \overline N$, and therefore Theorem 3 asserts that $Y = N$, thus
$\alpha$ is surjective. 
	\bigskip
{\bf Corollary.} {\it Let $\epsilon\: X \to N$ be an epimorphism and $N \subseteq Y$
an inclusion map with semisimple cokernel such that the
composition $X \to N \to Y$ is kernel-determined. 
Then there is an inclusion map $Y \to Z$ such that the composition $X \to Y \to Z$ 
has semisimple cokernel, is kernel-determined and satisfies}
$$
  |Z| = |N|+|\overline N|-|I_\epsilon(N)|.
$$  
In particular, the length of $Z$ only depends on $\epsilon.$
	\medskip
Proof: We can assume that $Y$ is a submodule of $\overline N$.  
Choose $N \subseteq Z \subseteq \overline N$ maximal with 
$Z \cap I_\epsilon(N) = N.$ According to Theorem 3, the composition
$X \to Y \to Z$ (which is the composition of $\epsilon$ and the inclusion map $N \to Z$)
is kernel-determined. The maximality of $Z$ implies that $Z + I_\epsilon(N) = 
\overline N$. The stated equality comes from the formula
$$
 |Z| + |\overline\epsilon(X)| = |Z\cap \overline\epsilon(X)| + |Z + \overline\epsilon(X)|.
$$
	\bigskip
{\bf Summery.} The kernel-determined morphisms can be characterized as suitable prolongations
of epimorphisms. Here, we call 
the composition $X \to Y \to Z$ a {\it prolongation} of $X \to Y$ provided
the map $Y \to Z$ is an inclusion map; 
the prolongation is said to be {\it proper} provided the
map $Y \to Z$ is a proper inclusion map.
	\medskip
(a) { Any epimorphism $X\to N$ is kernel-determined.}
	\medskip
(b) { If the map $X \to Y$ has a prolongation $X \to Y \to Y'$ 
which is kernel-determined, then $X \to Y$
is kernel-determined and $Y \to Y'$ is an essential extension.}
	\medskip
(c) { Let $X \to N$ be an epimorphism, and 
$N \subseteq Y \subseteq I(N)$. If $X \to N \to Y\cap \overline N$ is kernel-determined,
also $X \to N \to Y$ is kernel-determined.}
	\medskip
(d) { Any kernel-determined map $X \to Y$ has a 
maximal kernel-determined prolongation 
$X \to Y \to Y'$.}
	\medskip
(e) { If $X \to N$ is an epimorphism, and 
$N \subseteq Y \subseteq I(N)$, then $X \to N \to Y$ is kernel-determined if and only if
$Y\cap I_\epsilon(N) = N.$}
	\medskip
(f) { If $X \to N$ is an epimorphism and $X \to N \to Y$ is a 
maximal kernel-determined prolongation, then 
$$
  |\soc (Y/N)| =  |\soc (I(N)/N)| - |I_\epsilon(N)/N|;
$$
in particular, the length of $\soc (Y/N)$ is determined by $\epsilon$.}
	\medskip
Thus, { if $X \to N$ is an epimorphism
and $X \to N \to Y$ and
$X \to N \to Y'$ are  maximal kernel-determined prolongations,
then $\soc(Y/N)$ and $\soc(Y'/N)$ have the same length,} 
but $Y$ and $Y'$ may have different length, as the following example shows: 
	\medskip
{\bf Example 6.} Consider the representations of the following quiver with relations over the
field $k$: 
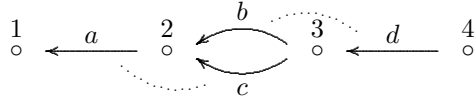
$$
{\beginpicture
\setcoordinatesystem units <2cm,1cm>
\multiput{$\circ$} at 0 0  1 0  2 0  3 0 /
\put{$1$} at 0 0.3
\put{$2$} at 1 0.3
\put{$3$} at 2 0.3
\put{$4$} at 3 0.3
\put{$a\strut$} at 0.5 0.2
\put{$b\strut$} at 1.5 .5
\put{$c\strut$} at 1.5 -.5
\put{$d\strut$} at 2.5 0.2
\arr{0.8 0}{0.2 0}
\arr{2.8 0}{2.2 0}

\arr{1.25 0.16}{1.2 0.1}
\setquadratic
\plot 1.2 0.1  1.5 0.3  1.8 0.1 /

\arr{1.25 -.16}{1.2 -.1}
\plot 1.2 -.1  1.5 -.3  1.8 -.1 /
\setdots <1mm>
\plot 0.7 -0.2  1 -.5  1.3 -.4 /
\plot 1.7   .4  2  .5  2.3  .2 /

\endpicture}
$$
We denote the simple, projective, or injective module corresponding to the vertex $x$
by $S(x), P(x), I(x)$, respectively. The full subquiver with vertices $2,3$ is the Kronecker
quiver, the representations with support in this subquiver will be said to be Kronecker modules. 
The 2-dimensional indecomposable Kronecker module which
is annihilated by $\lambda_1b+\lambda_2c$ (not both $\lambda_1,\lambda_2$ equal to zero) will
be denoted by $R(\lambda_1b+\lambda_2c).$  For example, $I(1)/S(1) = R(c)$ and
$\rad P(4) = R(b).$

Let $X = P(2)$ and $N = S(2)$ and $\epsilon\:X \to N$  the canonical projection $P(2)\to S(2)$. 
Then $\overline X = I(1)$ is indecomposable with composition factors $S(1),S(2),S(3)$.
The module $\overline N$ has length 3, namely one composition factor $S(2)$ and two composition factors $S(3)$, it is just the indecomposable injective Kronecker module of length $3$
and $I_\epsilon(N) = R(c).$
 
In view of Theorem 3, we are interested in the submodules $Y$ of $\overline N$ which
satisfy $Y \cap I_\epsilon(N) = N,$ thus $Y \cap R(c) = N.$
Besides $N$ itself, these are the Kronecker modules of the from $R(b+\lambda c)$
with $\lambda\in k.$ The modules $Z = R(b+\lambda c)$ provide the 
maximal kernel-determined prolongations $X \to Y \to Z$ of $X \to N$
inside $\overline N$.

Now only the map $X \to N \to R(b)$ has a proper kernel-determined 
prolongation, namely $X \to R(b) \to P(4)$.
The other maps $X \to N \to R(b+\lambda c)$ with  $\lambda\neq 0$ 
have no proper kernel-determined prolongation.
	\bigskip\bigskip
{\bf 6. References.}
	\medskip
\frenchspacing
\item{[A1]} Auslander, M.: Functors and morphisms determined by objects. In: 
Representation Theory of Algebras. Lecture Notes in Pure Appl. Math. 37.
  Marcel Dekker, New York (1978), 1-244. Also in: Selected Works of Maurice Auslander,
  Amer. Math. Soc. (1999).
\item{[A2]} Auslander, M.: Applications of morphisms determined by objects.
In: 
Representation Theory of Algebras. Lecture Notes in Pure Appl. Math. 37.
  Marcel Dekker, New York (1978), 245-327. Also in: Selected Works of Maurice Auslander,
  Amer. Math. Soc. (1999).
\item{[ARS]} Auslander, M., Reiten, I., Smal\o, S.:
  Representation Theory of Artin Algebras. Cambridge Studies in Advanced Mathematics 36.
Cambridge University Press. 1997.
\item{[K]} Krause, H.: Morphisms determined by objects in triangulated
  categories. arXiv:1110.5625.
	\bigskip\bigskip

\parindent=0truecm
Claus Michael Ringel \par
Department of Mathematics, Shanghai Jiao Tong University \par
Shanghai 200240, P. R. China, and  \par
King Abdulaziz University, \par
P O Box 80200, Jeddah, Saudi Arabia.\par

{\tt e-mail: ringel\@math.uni-bielefeld.de}
\bye